\documentclass{article}

\usepackage{latexsym}
\usepackage{epsfig}  
\begin{document}

\title{\huge A Note on Curve Counting Scheme in an Algebraic Family and
 The Admissible Decomposition Classes}
\author{Ai-Ko Liu\footnote
{email address: akliu@math.berkeley.edu, Current Address: 
Mathematics Department of U.C. Berkeley}}
\maketitle

\newtheorem{theo}{Theorem}
\newtheorem{lemm}{Lemma}
\newtheorem{prop}{Proposition}
\newtheorem{rem}{Remark}
\newtheorem{cor}{Corollary}
\newtheorem{mem}{Examples}
\newtheorem{defin}{Definition}
\newtheorem{axiom}{Axiom}
\newtheorem{conj}{Conjecture}
\newtheorem{exam}{Example}
\newtheorem{assum}{Assumption}

\bigskip

 In this paper, we discuss the generalized scheme for curve
 counting in the family Seiberg-Witten theory.
 Even though the original motivation is to study the
 Mcduff's proposal in $b^+_2=1$ category of symplectic four manifolds,
 we will formulate our scheme in an algebraic(Kahler) set up. The
 material considered in this paper will be relatively
 elementary. Nevertheless, the theory discussed here has played an essential
 role in the long paper [Liu1].

 As a major application of the discussion, one may apply our scheme
 in the proof of the G${\ddot o}$ttsche's conjecture about counting of
 holomorphic curves[Liu1]. The current 
scheme is motivated from the discussion [Mc] of
 pseudo-holomorphic curves in symplectic four manifolds. However we
 will restrict our discussion to the algebraic varieties here.
 One can translate our scheme to pseudo-holomorphic category
 by replacing holomorphic curves to pseudo-holomorphic curves.
 
\medskip

 Given an algebraic family of algebraic surfaces ${\cal X}\mapsto B$,
 the family Seiberg-Witten invariant ${\cal AFSW}$ (or $FSW$ in the smooth
 category) enumerates the algebraic curves (or pseudo-holomorphic curves
 in the symplectic category) within the family dual to a given cohomology class.
 The basic phenomenon we will study is that not only smooth curves may
 appear in the enumeration, curves contain multiple coverings of the so-called
 exceptional curves may also occur.  The general question we are interested
 at is

\noindent{\bf Question}: How to relate the contribution from the smooth 
 curves to the original family invariant? 

\medskip

 The general strategy to answer the question is to subtract the 
 contributions from the various configurations containing multiple coverings of
 exceptional curves. The purpose of the
 current paper is to provide a skeleton of the curve counting scheme.

\section{\bf The family scheme of Algebraic Surfaces}\label{section; AS}

\bigskip

In this section, consider ${\cal X}\mapsto B$ to be a relatively
 smooth algebraic fibration over a smooth base. 
For simplicity, the field of definition
can be taken to be ${\bf C}$. The same scheme will work for any
algebraic closed field of characteristic zero as well. Even though
 we do not aim at general symplectic four manifolds, we will
 recall the phenomena in the symplectic set up from time to time.

 The fibers of the fibration are taken to be smooth projective
 surfaces. In the following, we denote $dim_{\bf C}B$ to be the complex
 dimension of the base even though the scheme work well even when
 $B$ does not carry complex structures.

 First we introduce certain notations. 
Let $C$ denote the cohomology class represented by holomorphic curves.
Let $d_{GT}(C)={C^2-c_1(K_X)\cdot C\over 2}$
 denote the complex Gromov Taubes dimension [T3] in Taubes theory.
 Then the family Gromov-Taubes dimension
 of $C$ is defined axiomatically to be $d_B(C)=d_{GT}(C)+dim_{\bf C} B$.

 Given the class $C$ as the sum of $C_1$ and $C_2$, one has the
 following equality between their formal family Gromov-Taubes dimensions.

$$d_B(C)=d_B(C_1)+d_B(C_2)-dim_{\bf C}B.$$

 This equality reflects that the family moduli space associated to
 $C$ is the fiber product of those of $C_1$ and $C_2$.

 One can easily generalize the equality to
 more than two $C_i$ and the new equality is 
 
$$d_B(\sum_{i\in I} C_i)=\sum_{i\in I} d_B(C_i)-(|I|-1)dim_{\bf C}B,$$
with $|I|$ being the cardinality of the index set $I$.

\medskip

\begin{defin}\label{defin; exception}
 A cohomology class $e\in H^2(M, {\bf Z})$ is said to be an exceptional class
 if 

 (i). $e$ is a primitive element in the lattice $H^2(M, {\bf Z})$.

 (ii). $e^2=e\cup e[M]<0$.

 A pseudo-holomorphic curve in an almost complex four-manifold 
$M$ is said to be an exceptional curve if
 it is poincare dual to an exceptional class $e$.
\end{defin}

 We have the following proposition regarding irreducible exceptional curves,

\begin{prop} \label{prop; unique}
 Let $\Sigma_1, \Sigma_2$ be two Riemann surfaces and let $f_1:\Sigma_1\mapsto M$
 $f_2:\Sigma_2\mapsto M$ be two (pseudo-)holomorphic maps into the almost complex
 four-manifold $M$ with $(f_i)_{\ast}[\Sigma_i]$ dual to an exceptional class
 $e\in H^2(M, {\bf Z})$, then $\Sigma_1=\Sigma_2$ and 
$f_1$ and $f_2$ coincide.
   Namely, $f_i(\Sigma_i), i=1, 2$ coincide and is the only irreducible
 (pseudo-)holomorphic curve dual to $e$ in $M$.
\end{prop}

\medskip

\noindent Proof of the proposition:  Because $e$ is primitive, the maps 
 $f_i$ from $\Sigma_i$ to $f_i(\Sigma_i)$ are of degree $1$. From [Mc],
 $f_i(\Sigma_i), i=1, 2$ have 
at most a finite number of isolated singularities and  
 $f_i$ are immersions away from the singularities in the images.
 Because $\Sigma_i$ are irreducible, $f_i(\Sigma_i)$ are irreducible, too.

 We argue that $f_1(\Sigma_1)=f_2(\Sigma_2)$.  If not, the
 sets $f_1(\Sigma_1)=\cap f_2(\Sigma_2)$ is of a finite cardinality.
 However, again by [Mc2] each intersection point contributes positively
 to the total intersection number, the total intersection number should 
 be non-negative. On the other hand, $f_1(\Sigma_1)\cap f_2(\Sigma_2)
 =(f_1)_{\ast}[\Sigma_1]\cap (f_2)_{\ast}[\Sigma_2]=e\cup e[M]<0$.
 
 This gives us the necessary contradiction. Thus $f_1(\Sigma_1)=f_2(\Sigma_2)$.
 Once we know $f_1(\Sigma_1)=f_2(\Sigma_2)$ and both pseudo-holomorphic 
 maps are of degree one,
 Both $\Sigma_1, \Sigma_2$ can be re-constructed as the normalization of the
 complex curve $f_1(\Sigma_1)=f_2(\Sigma_2)$. Therefore $\Sigma_1=\Sigma_2$ and
 $f_1=f_2$. $\Box$

 \medskip

  Even though the ``irreducible'' curve dual to $e$ is always unique, there can be
 two or more reducible pseudo-holomorphic curves dual to $e$
 with more than one irreducible component. In this case, the conclusion is
 weaker and one
can only deduce that two curves share at least one irreducible component
 and the fundamental class of this irreducible component has a negative 
self-intersection number in the four-manifold $M$.

 Recall that Taubes theory [T3] asserts the equivalence of Seiberg-Witten
 invariant
 and a version of Gromov invariant for symplectic four
 manifolds. It indicates that the diffeomorphism invariants $SW$ is
 equivalent to the symplectic invariant $Gr$. Despite of the simplicity of
 the statement, the actual proof [T1], [T2], [T3]
 involves sophisticated analysis and 
 an amount of new ideas.

  It is less well known that the equivalence of $SW$ and $Gr$ fails
  for the general symplectic four manifolds with $b^+_2=1$.
  Originally Taubes asserted his theorem for $b^+_2>1$ category.
  Later he extended his theorem to $b^+_2=1$ case with some additional
 assumption. In the mean time, it was discovered experimentally by
 the current author in [LL] that the assertion would not be true
 without the additional assumption.  This motivated Mcduff to
 change the original definition of the Gromov invariant in order to match up
  with $SW$.

\medskip

\begin{defin}\label{defin; Taubes'}
 In defining the Seiberg-Witten invariants of $b^+_2=1$ symplectic 
 four-manifold $M$, the space of generic Riemannian metrics $g$ and self-dual
 two forms $\mu$ on $M$ are divided into chambers. Given a symplectic two form
 $\omega$ on $M$ and a Riemannian metric $g$ such that $\omega$ is self-dual
 with respect to $g$, the $(g, r\omega=\mu)$, $r\mapsto \infty$ determines
 a unique chamber, called the Taubes' chamber in the following discussion.
\end{defin}

\medskip

 Let us recall Taubes theorem for $b^+_2=1$ symplectic manifolds.
 
\begin{theo}(Taubes) \label{theo; =1}
 Let $C$ be a cohomology class in $H^2(M, {\bf Z})$, with $M$ being
 a symplectic four manifold with $b^+_2=1$.  Assume additionally
 that $C\cdot S\geq -1$ for all spherical class $S$ with $S^2=-1$.
 Then the statement $SW(2C-K_M)=Gr(C)$ holds for $C$, where $SW(2C-K_M)$ is
 evaluated in the Taubes' chamber.
\end{theo}
  
 By a spherical class $S$ with $S^2=-1$ one means that the class 
is represented by a $S^2$ with self intersection number $-1$. 

In the following, we give a simple example that the theorem does not
 hold without modification for the classes violating this condition.

\begin{mem}
 Consider $M$ to be ${\bf CP}^2\sharp \overline{{\bf CP}^2}$, the symplectic
 four manifold constructed by ${\bf CP}^2$ by blowing up one point in 
 ${\bf CP}^2$. 
 Let $H\in H^2(M, {\bf Z})$
 denote the (pull-back of) 
 the hyperplane class and let $E\in H^2(M, {\bf Z})$ denote the exceptional class.

 Take $C=3H+2E$,and a simple calculation shows that $d_{GT}(C)=9-1=8$.
 Through the calculation of wall crossing formula it is easy to see
 that $SW(9H+3E)=\pm 1$. On the other hand, we argue that 
a reasonable definition of
 Gromov invariant would be $Gr(3H+2E)=0$. The class 
 $3H+2E$ can not be represented by irreducible curves. The
 representatives are the disjoint union of a cubic curve dual to $3H$ 
along with the double covering of the exceptional curve $E$.
 Even though curves dual to $3H$ gives rise to nonzero invariant, the multiple
 covering of $E$ has a negative Gromov dimension $d_{GT}(2E)=-1$. Thus the
 total Gromov-Taubes invariant should be $Gr(3H)\times 0=0$.
\end{mem}

\medskip

 This simple example illustrates the subtlety to Taubes' theory.
 Even though the $SW(9H+3E)\not=0$ in the Taubes' chamber and the analysis in 
$SW\mapsto Gr$ [T1] still implies the existence of 
 pseudo-holomorphic curves in $C$, numerically the curves are counted zero
 in the Gromov-Taubes theory.  As similar type of phenomena appears for
 any non-minimal symplectic four manifolds with $b^+_2=1$, it becomes
 the major topological obstruction to identify $SW$ and $Gr$.
 The way Taubes dealed with this problem is to rule out the classes $C$
 which potentially can be represented as disjoint unions of
 pseudo-holomorphic curves and some multiple coverings of exceptional $-1$
 curves. This explains the extra condition $C\cdot S\geq -1$ in Taubes'
 theorem.

 Even though it is not completely obvious. this ill symptom is closely
 related to the fact that $b^+_2=1$ symplectic manifolds are not
 of simple type in Taubes' chamber. In the fundamental paper of
 Taubes, symplectic four manifolds with $b^+_2>1$ were  proved to be
 of simple type. Thus, the ill symptom does not occur to them. 

 On the other hand, it is a consequence of the wall crossing formula
 [LL] that the non-simple type of $b^+_2=1$ symplectic manifolds is
 directly related to the non-vanishing of the wall crossing numbers.
  
 The primitive goal to develop the family Seiberg-Witten theory is to
 discuss the family Seiberg-Witten theory in Taubes' chambers (defined by
 a large perturbation of fiberwise self-dual symplectic forms).
 As a similar application of the family wall crossing formula implies
 the non-simple type-ness, one would expect that the similar failure
 of $SW=Gr$ would occur. Unlike the $B=pt$ case that $-1$ curves are
 the only pseudo-holomorphic curves which persist, the topological
 types of the exceptional
 curves which persist in the family are less restricted.

 Viewed from a different angle, the major distinction between the
 standard Gromov-Witten theory and Gromov-Taubes theory lies in the
 fact that the latter theory does not restrict the topological types
 of the pseudo-holomorphic maps.

 The combination of these two issues make the identification between
 $FSW$ and $FGr$ extremely difficult. 
I.e./ within a given family of symplectic manifolds and a fiberwise 
 monodromy invariant class $C$, 
there can be a whole 'zoo' of exceptional 
 curves which may appear in some pseudo-holomorphic curve 
representations of the class $C$. To make sense of Gromov-Taubes
  invariants, one has to deal with these exotic objects in a more
 systematical way.

\bigskip

\begin{mem} \label{mem; explain}
  For the readers with a background of Gromov-Witten theory, it should be
 cautious not to think of Gromov-Taubes invariants as identical to the
 standard Gromov-Witten invariants. Besides the question of allowing 
dis-connected domain curves, the expected Gromov-Taubes 
dimensions of multiple covering
 of exceptional curves are different from the dimension formulae of 
Gromov-Witten invariants viewed as multiple covering pseudo-holomorphic maps
into the symplectic manifold.

  For simplicity, let $e\in H^2(M, {\bf Z})$
 be an exceptional class $e^2=-k, c_1(K_M)\cdot e=-2+k$ representing an
 exceptional sphere. Let $m\in {\bf N}$
 be a positive integer, the expected Gromov-Taubes dimension of the class
 $me$ is given by $(-m^2-m)k+2m$. 

 On the other hand, the expected dimension of the pseudo-holomorphic maps
 from $S^2$ to $M$, dual to $me$, (modulo diffeomorphisms on $S^2$) is given by 

  $$2\{c_1(M)\cdot C-(g-1)\cdot dim_{\bf C}M+3g-3\}=2\{-c_1(K_M)\cdot C+2-3\}
=-2mk+2m-2.$$
 
 The former is quadratic with respect to $m$, while the usual Gromov-Witten
 expected dimension is linear in $m$.

 Thus, even though we will still call pseudo-holomorphic curves dual to $me$ a
 multiple covering of exceptional curves in $e$, we advise the readers not
to confuse them with the multiple covering in the sense of maps. Instead, it
 is wiser to think of it as $m$ copies of pseudo-holomorphic curves dual to $e$
 sitting on top of each other.
\end{mem}

\medskip

\section{\bf The Pointwise Calculation of Family Dimension}
\label{section; point}

\medskip

 Firstly, let us review the original dimension count argument of Taubes.
 
 Let $(M, \omega)$ be a symplectic four-manifold with a compatible almost
 complex structure and $C\in H^2(M, {\bf Z})$ be a cohomology class with
 a positive energy $C\cdot \omega>0$. In Taubes' theory, he allows the
 pseudo-holomorphic curves to have more than one irreducible component.
 
Suppose that there is a pseudo-holomorphic curve poincare dual to $C$.
 Then there is a finite collection of Riemann surfaces $\Sigma_i, 1\leq i\leq k$
 and the pseudo-holomorphic maps $f_i:\Sigma_i\mapsto M$ such that
 $\sum_{i\leq k} (f_i)_{\ast}[\Sigma_i]\in H_2(M, {\bf Z})$ is poincare dual
 to $C$.

  In case the map $f_i$ is of degree $m_i$, we may write 
$PD((f_i)_{\ast}[\Sigma_i])$ as $m_ie_i$, $m_i\in {\bf N}$. Then we have
 the following equality 

 $$C=\sum_{i\leq k} m_i e_i.$$

 Conversely, if $C$ is written as $\sum_{i\leq k} m_i e_i$ and each of
 $e_i$ is represented by an irreducible pseudo-holomorphic curve, then one 
 takes the union of them (counting multiplicity) and represent $C$ as a
 pseudo-holomorphic curve in $M$.

 Taubes would like to study all the possible decompositions of $C$ into
 different $e_i$ which will survive under generic compatible 
 almost complex structures
 perturbation of $M$.

  One makes three additional assumptions on $e_i$,

\medskip

 (i). $e_i\cdot \omega>{\cal E}(M, \omega)>0$ for some manifold dependent
 lower bound of harmonic energy.

\medskip

 (ii). The expected Gromov-Taubes dimension 
$d_{GT}(e_i)={e_i^2-e_i\cdot c_1({\bf K}_M)\cdot e\over 2}\geq 0$.

\medskip

(iii). $e_i\cdot e_j\geq 0$ for all $i\not= j$.

\medskip

(iv).  $e_i^2+c_1({\bf K}_M)\cdot e_i=2g_{arith}(e_i)-2\geq -2$.

\medskip

\begin{defin}\label{defin; decomposition}
 Let $e_i, i\leq k$ be a finite number of classes in $H^2(M, {\bf Z})$ 
 which satisfy $(i)., (ii)., (iii).$ and $(iv).$
  The expression $C=\sum_{i\leq k} m_i e_i, m_i\in {\bf N}$ is called 
 a (cohomological) decomposition of $C$. 

\end{defin}

\medskip

 The reason that one imposes $(i)$ is because a pseudo-holomorphic curve dual to 
 $e_i$ always has a positive energy.

 If an irreducible irreducible 
curve dual to $e_i$ has negative Gromov-Taubes dimension 
$d_{GT}={e_i^2-c_1({\bf K}_M)\cdot e_i\over 2}<0$, then by Fredholm theory
 this type of curves may disappear after a generic perturbation of compatible 
 almost complex structures on $M$, which does not have chance to 
contribute to the Gromov invariant defined by Taubes [T3].

 From [Mc2], one knows that two distinct irreducible pseudo-holomorphic curves
 in an almost complex four-manifold $M$ 
intersect positively. Thus $e_i\cdot e_j$, the sum of all the local
 intersection contribution, should be non-negative as well.

 Because $e_i$ is represented by an irreducible pseudo-holomorphic curve on
 an almost complex four-manifold $M$, it satisfies the adjunction formula
 with $g_{arith}$ being the arithmetic genus of the curve.

 Let $C$ be written as 
 $\sum_{i\leq k} m_i e_i$ with multiplicity
 $m_i\geq 1$ satisfying $(i).$, $(ii).$,  $(iii).$ and $(iv)$.

Then $$2d_{GT}(C)={C^2-C\cdot c_1({\bf K}_M)}=
\bigl(\sum_{i\leq k} (e_i^2-e_i\cdot c_1({\bf K}_M))$$
$$+2\sum_{i\not= j} m_i m_j e_i\cdot e_j+
\sum_{i\leq k}((m^2_i-1)e_i^2+(1-m_i)e_i\cdot c_1({\bf K}_M))\bigr)$$

$$=2\sum_{i\leq k}d_{GT}(e_i)+2\sum_{i\not= j} m_i m_j e_i\cdot e_j+
\sum_{i\leq k}\bigl((m^2_i-m_i)e_i^2+
(m_i-1)(e_i^2-e_i\cdot c_1({\bf K}_M))\bigr).$$

 By the assumption (iv). $e_i^2+c_1({\bf K}_M)\cdot e_i\geq -2$ and 
 by $(ii).$ $2d_{GT}(e_i)=e_i^2-c_1({\bf K}_M)\cdot
 e_i\geq 0$. Then we know that $e_i^2\geq -1$. Suppose $e_i^2=-1$. From
 $(ii).$ again we get $-1=e_i^2\geq c_1({\bf K}_M)\cdot e_i$. Then
 $-2\geq e_i^2+c_1({\bf K}_M)\cdot e_i$ and $c_1({\bf K}_M)\cdot e_i=-1$
 as well. If it is the case, $g_{arith}(e_i)=0$ and one can argue that
 the pseudo-holomorphic curve representing $e_i$ must be a so-called $-1$ curve.
  
 Because Taubes' goal is to develop a version of Gromov invariant 
 which can be identified with $SW(2C-c_1({\bf K}_M))$ (the $spin^c$ class
 $2C-c_1({\bf K}_M)$ is in an additive notation), he is able to
 use the Seiberg-Witten simple type-ness condition on $b^+_2>1$ symplectic
 four-manifolds and the blowup formula of Seiberg-Witten invariants 
 to deduce $m_i=1$ for all such $-1$ classes $e_i^2=-1$.

 Thus, one find that the last term in the expansion of $2d_{GT}(C)$ is
 always non-negative.

 In order to count pseudo-holomorphic curves dual to $C$, one imposes
 $d_{GT}(C)$ number of generic points and require the pseudo-holomorphic
 curves dual to C to pass through these generic points.

 From above we find that $$d_{GT}(C)-\sum_{i\leq k}d_{GT}(e_i)\geq 0.$$

  If the difference $d_{GT}(C)-\sum_{i\leq k}d_{GT}(e_i)$ 
is strictly positive, then by dimension reason 
there can be no pseudo-holomorphic 
curves in $\sum_{i\leq k} e_i$ which pass through all these $d_{GT}(C)$ points, 
 after one adopts the Fredholm argument to perturb the almost complex structures
 of $M$. In order $\sum m_i e_i$ contributes to the Gromov invariant, the
 non-negative
 sum $2\sum_{i\not= j} m_i m_j e_i\cdot e_j+
\sum_{i\leq k}\bigl((m^2_i-m_i)e_i^2+
(m_i-1)(e_i^2-e_i\cdot c_1({\bf K}_M))\bigr)$ has to vanish term by term. 

 That is to say, 

$(a).$ $m_i=1, \forall i\leq k$.

$(b).$ $e_i\cdot e_j=0, \forall i\not=j$.

$(c).$ Each curve dual to $e_i$ must be smooth.

\medskip

 In other words, distinct $e_i$ cannot intersect. Each irreducible 
pseudo-holomorphic curve
 appears with multiplicity one.

 One may develop the following 'philosophical' idea which helps to 
 explain what happens.

\medskip

(a)'.  If $m_i>1$ for some $e_i^2>0$, ideally one may choose two
 distinct pseudo-holomorphic curves dual to $e_i$ and they intersect
 positively. Then the smoothing of all these
 intersection singular points produces an irreducible curve dual to $2e_i$, whose
 dimension $d_{GT}(2e_i)>d_{GT}(e_i)+d_{GT}(e_i)$.  Continue in this fashion,
 Seiberg-Witten theory is expected to count irreducible curves in $m_ie_i$
 rather than $m_i$ copies of curves dual to $e_i$, which formally can be viewed as
 a degeneration from an irreducible multiplicity one curve dual to $m_ie_i$.

\medskip

(b)'. If $e_i\cdot e_j>0$ for some $i\not j$,
 one may think of the smoothing of the $e_i\cdot e_j$ 
intersection points (counted with multiplicity) and consider (formally) a
 curve dual to $e_i+e_j$. Then the union of the curves dual to $e_i$ and $e_j$
 can be thought as a degeneration of some irreducible curve dual to $e_i+e_j$
 and $d_{GT}(e_i+e_j)>d_{GT}(e_i)+d_{GT}(e_j)$.

\medskip

(c)'. If a curve dual to $e_i$ develops certain singularities, it can
 be thought of a degeneration of the smooth curves dual to $e_i$ satisfying
 the adjunction equality $e_i^2+c_1({\bf K}_M)\cdot e_i=2g(e_i)-2$. The
 curves with singularities are of lower expected dimension than the expected 
 Gromov-Taubes dimension $d_{GT}(e_i)$.

\medskip

In the $b^+_2=1$ category, $m_i\geq 1$ for the $-1$ classes $e_i^2=-1$.
 Then the same argument breaks down.

One way to remedy is to impose extra condition on $C$, as was done by
 Taubes (see the statement of theorem. \ref{theo; =1}).

Mcduff introduces a different way to remedy the situation. She (see [Mc])
 has shown that

\begin{prop}(Mcduff) \label{prop; mcduff}
 Let $M$ be a $b^+_2=1$ symplectic four-manifold and let $C\in H^2(M, {\bf Z})$
 be a class satisfying $d_{GT}(C)={C^2-C\cdot c_1({\bf K}_M)\over 2}\geq 0$,
 $C\cdot \omega>0$.  Then there exists a finite number of $-1$ classes,
 $e_i$, $e_i^2=-1$ satisfying the following conditions:

(i). Each $e_i$ is represented by a $-1$ pseudo-holomorphic curve.

\medskip

(ii). $e_i\cdot e_j=0$ for $i\not=j$.

\medskip

(iii). $C\cdot e_i=-n_i<0$.

 Then one may re-write $C=(C-\sum n_ie_i)+\sum n_i e_i$ and $C-\sum n_ie_i=C_{red}$
 is perpendicular to all these $e_i$, namely $C_{red}\cdot e_i=0$.

\end{prop}

\medskip

 Mcduff proposed [Mc] to define $Gr(C)$ using the class $C_{red}$ instead of $C$.

\medskip

\begin{defin}\label{defin; taubes}
 Let $C=\sum_i m_i e_i$ be a cohomological decomposition of $C$. The
 decomposition is said to be of Taubes' type if 

(a). $e_i\cdot e_j=0$ for $i\not j$.

\medskip
 
(b). $m_i=1$ for all $i$.

\end{defin}

\medskip

In the following, we generalize the Mcduff's proposal to the
 family case. It turns out all the different possibilities of 
decompositions of curves appearing in the
 $B=pt$ case can also appear in the family case. Moreover, in the family
 case there are many new possible decompositions of curve classes
 which are absent in 
the $B=pt$ case due to dimension reason.

  Let $\pi:{\cal X}\mapsto B$ be a fiber bundle of symplectic four-manifolds with 
 the relative symplectic form $\omega_{{\cal X}/B}$. 

 Similar to the condition $(i). (ii). (iii).$ and $(iv).$ for the $B=pt$ cases,
 one imposes the following conditions on the classes $e_i$.

\medskip

(Fi).$e_i\cdot \omega_{{\cal X}/B}>{\cal E}({\cal X}, \omega_{{\cal X}/B})>0$ for 
some fiber bundle dependent lower bound of harmonic energy.

\medskip

(Fii).The expected family Gromov-Taubes dimension 
 $dim_{\bf C}B$+$d_{GT}(e_i)={e_i^2-e_i\cdot c_1({\bf K}_M)\cdot e\over 2}\geq 0$.

\medskip

(Fiii). $e_i\cdot e_j\geq 0$ for all $i\not=j$.

\medskip

(Fiv). $e_i^2+c_1({\bf K}_{{\cal X}/B})\cdot e_i=2g_{arith}(e_i)-2\geq -2$.

\medskip

At this moment we are doing the pointwise analysis for different $b\in B$,
 we do not take into account
 the monodromy action of $\pi_1(B, b_0)\mapsto H^2(\pi^{-1}(b_0), 
 {\bf Z})$. At times (e.g. for the universal families $M_{l+1}\mapsto M_l$),
 the monodromy representation is completely trivial and we can ignore 
 it. Otherwise, we have to consider the equivalent classes of $e_i$ 
or decompositions under the action of 
$\pi_1(B, b_0)\mapsto H^2(\pi^{-1}(b_0), {\bf Z})$.

  Let us discuss how does the ordinary Taubes' dimension count argument
 generalized to the family case. Because we will use the relative 
 canonical bundle ${\bf K}_{{\cal X}/B}$ throughout the discussion, we
 will skip the subscript ${\cal X}/B$ and denote it by ${\bf K}$.

Using the family dimension formula,

$$ dim_{\bf R}B+2d(C)=dim_{\bf R}B+{C^2-C\cdot c_1({\bf K})}
=(\sum_{i\leq k} (dim_{\bf R}B+e_i^2-e_i\cdot c_1({\bf K}))-(k-1)dim_{\bf R}B$$
$$+2 \sum_{i\not= j}
 m_i m_j e_i e_j+\sum_i((m^2_i-1)e_i^2+(1-m_i)e_i\cdot K)).$$
\medskip

  In a given family, the condition $(Fii).$ is weaker than the condition 
 $(ii).$ at the $B=pt$ case. Thus, there may be some $e_i$ with $e^2_i<-1$.
 Then the last expression $\sum_i((m^2_i-1)e_i^2+(1-m_i)e_i\cdot K))$ may be
 negative for $m_i\not=1$.  In other words, the formal dimension
 expected from family Seiberg-Witten theory (equal to $dim_{\bf R}B+
 2d_{GT}(C)$ can be smaller than the actual
 formal dimension on the Gromov-Taubes side, $(\sum_{i\leq k} 
(dim_{\bf R}B+e_i^2-e_i\cdot c_1({\bf K}))-(k-1)dim_{\bf R}B$.

 \bigskip

Motivated from Mcduff's proposal [Mc] and Taubes' theorem \ref{theo; =1}, 
 let us denote $P$ as the index subset $P\subset \{1, 2, \cdots, k\}$ 
such that $e_i\cdot C<0, i\in
 P$, with 
 $e_i^2<0$. Then we can always regroup the decomposition of $C$ as 

 $$C=F+E,F=\sum_{i\notin P}m_i e_i; E=\sum_{j\in P}m_j e_j.$$

  Namely, one can view $F$ as a whole without going into
the details of the decomposition of the class $F$.

 Then the previous expression can be expanded easily into the following 

$$dim_{\bf R}B+2d_{GT}(C)=F^2-c_1({\bf K})\cdot F+2F\cdot E+E^2-c_1({\bf K})
\cdot E+dim_{\bf B} B.$$

 Then we rewrite the  term $F\cdot E$ term into $F \cdot
 E=(C-E)\cdot E$, then we have

$$dim_{\bf R}B+2d_{GT}(C)
=F^2-c_1({\bf K})\cdot F+2(C-E)\cdot E+E^2-c_1({\bf K})\cdot E+dim_{\bf R} B.$$

 Let us collect $F^2-c_1({\bf K})\cdot F+dim_{\bf R}B+
 \sum_{i\in P} (e_i^2-c_1({\bf K})\cdot e_i)$ into
 a single term. It is the family dimension of $F$, along with
 the all the $e_i, i\in P$.
The sum of the left-over terms, called the family dimension discrepancy, 
 has the following form

$$\Delta_C(E)=2C\cdot E-E^2-c_1({\bf K})\cdot E-
\sum_{i\in P}(e_i^2-c_1({\bf K})\cdot e_i).$$

 To study how do the multiplicities $m_i, i\in P$ in $E=\sum m_i E_i$ 
affect the family dimension, we have to introduce some combinatorial language. 
 combinatorial language. 

\medskip

\begin{defin} \label{defin; exc}
 Let $C$ be a monodromy invariant fiberwise cohomology class of
 $\pi: {\cal X}\mapsto B$. Given a point $b\in B$, 
Let $e_i, i\in P$ be all the classes represented by irreducible 
 pseudo-holomorphic curves over $b$ with $C \cdot e_i<0, e_i^2<0$. Then define
 the exceptional cone of $C$ over $b$ to be 
 ${\cal EC}_{b}(C)=\sum_i{\bf R}^{\geq 0} e_i\in H^2(\pi^{-1}(b), {\bf R})$.
\end{defin}

\medskip

As $E$ is written as$\sum m_i e_i$ with $m_i\geq 0$. We can view
 E as an element in the cone ${\cal EC}_b(C)$
  generated by these $e_i$. 

  Given such a cone ${\cal EC}_{b}(C) \in H^2(\pi^{-1}(b),{\bf R})$ 
with an indefinite intersection
 form, it is possible to define the dual cone ${\cal EC}_b^{*}(C)$ to be
 the elements in $H^2(X,{\bf R})$ which have non-negative intersection 
pairings with elements
 in ${\cal EC}_b(C)$. As ${\cal EC}_b(C)$ is usually not a top dimensional
 cone in $H^2$, 
 the dual cone ${\cal EC}_b^{\ast}(C)$ usually is the direct sum of a vector space
 cone with a reduced dual cone in the minimal subspace containing ${\cal EC}_b(C)$.

 As a preparation, we want to prove some simple lemma characterizing the cone
 ${\cal EC}_b(C)$.

 Let us review some definitions.

\begin{defin} \label{defin; adm}
 Let ${\cal EC}_b(C)$ be an exceptional cone as
 described before. Then ${\cal EC}_b(C)$ is said to be admissible over $b$
 if $(C-{\cal EC}_b(C))\cap {\cal EC}_b(C)^{\ast}$ contains at least
 one lattice point.
\end{defin}

\medskip

 The following proposition clarifies the relation between the intersection
 form on $H^2(\pi^{-1}(b), {\bf R})$ with the admissible cone ${\cal EC}_b(C)$.

\begin{prop}\label{prop; negative}
 Suppose that the exceptional cone over $b$, ${\cal EC}_b(C)$, is
 admissible. Then the restriction of the 
 intersection quadratic form on the cone ${\cal EC}_b(C)$ is negative definite.
\end{prop}

\medskip

 This proposition implies that the term $-E^2$ in $\Delta_C(E)$ always 
contributes positively when $E\in {\cal EC}_b(C)$.

 Notice that we do not claim that the intersection form is negative 
definite on the whole minimal vector space in $H^2$
 containing ${\cal EC}_b(C)$. If the fibration $\pi:{\cal X}\mapsto B$
 is algebraic. Then all the $e_i$ are of type $(1, 1)$.
 Recall that Hodge index theorem asserts the intersection form has only one
 positive eigenvalue. In this case
this proposition asserts that the exceptional cone
 is disjoint to the forward and backward light cones.

\medskip

 The proposition is a generalization of Mcduff's proposal.

\medskip

\begin{rem}\label{rem; mc}
 When $B=pt$, the only exceptional curves satisfying $(ii).$ are $-1$ curves.
 Suppose $e_i, e_j\in {\cal EC}_{pt}(C)$ are two different $-1$ curves. 
 Then the proposition implies that $e_i+e_j\in {\cal EC}_{pt}(C)$ is
 of negative square. In other words, $(e_i+e_j)^2=-2+2e_ie_j<0$. Then $e_i\cdot 
 e_j$ must be $0$. This recovers Mcduff's proposal that all $-1$ curves
 $e_i, e_i\cdot C<0$ must be perpendicular to each other.
\end{rem}

\medskip

\noindent Proof of the Proposition: Set $|P|=n$. 
 As ${\cal EC}_b(C)$ is admissible, then there must be some tuple 
 $(m_1,m_2,\cdots, m_n)\in {\bf N}^{n}$ 
  such that $(C-\sum_{1\leq i\leq n}m_i e_i)\cdot e_j \geq 0$ for all 
 $j\in P$. 
 For simplicity let us re-scale and use the $\underline{e}_i=
m_i e_i$ as the new generators.
  Notice that it is related to the old one by a positive scaling. 
  Then we have the following inequality

 $$(\underline{e}_1+\underline{e}_2+\cdots +\underline{e}_n)
\cdot \underline{e}_i\leq C\cdot \underline{e}_i<0,$$
 for all $i\in P$. Now we use that
 $\underline{e}_i\cdot
 \underline{e}_j\geq 0$
 for $i\not= j$. 

Then we must have more inequalities of the similar
 type. Let $S\subset P$ be any nonempty subset of the index set $P$.
 Then $$(\sum_{i\in S} \underline{e}_i)\cdot \underline{e}_j
<-(\sum_{i\in P-S}\underline{e_i}\cdot \underline{e_j})<0, j\in S.$$

\begin{lemm}
Given any two index subsets $A, B\subset P$ such that one includes the
other. Then it follows that $\{\sum_{i\in A}\underline{e}_i\} \cdot
 \{\sum_{j\in B}\underline{e}_j\}<0$.
\end{lemm}
\medskip  

\noindent Proof of the Lemma: By symmetry we may assume $A\supset B$. Then 
 we take $A=S$ in the above discussion and we have 

 $$(\sum_{i\in A} \underline{e}_i)\cdot \underline{e}_j<0, j\in A.$$

 We may choose $j\in B\subset A$ and sum up all these inequalities for
 $j$ running through $B$, we get

 $$(\sum_{i\in A} \underline{e}_i)\cdot (\sum_{j\in B}\underline{e}_j)<0.$$ $\Box$

 Now we are ready to prove the statement in the proposition.
 Let $x$ be any element in the cone ${\cal EC}_b(C)$. Then it can be written as
 $\sum_{i\in P} c_i \underline{e}_i$ with $c_i\in {\bf R}+$.

  For convenience let us rearrange the indexes in $P$ in such a way 
that the coefficients $c_i$ are monotonically decreasing for increasing $i$.
  After such a permutation of indexes,
  let us consider elements of the form $f_a=\sum_{a\geq j\geq 1} 
\underline{e}_j$.
 Then the same element $x$ can be written alternatively as 

 $$E=\sum_{l\leq n}(c_l-c_{l+1})f_l,$$
 where we have set $c_{n+1}=0$. Therefore $E$ is written as an effective
 (because $c_l-c_{l+1}\geq 0$ expression over the new generators $f_ls$.

  On the other hand, the index sets involve in defining $f_l$ form a
 monotonic chain. They are of the form
 $\{1\}\subset \{1, 2\}\subset \{1, 2, 3\}\cdots $.

Therefore $f_l\cdot f_{l'}<0$ for all pairs of 
 $(l,l')$.
  Then it is easy to see that 
$$E\cdot E=\sum_{l, l'\leq n}(c_l-c_{l+1})(c_{l'}-c_{l'+1})f_l\cdot f_{l'}<0.$$

 The equality
 can hold only when $c_l-c_{l+1}=0$ for all $l$. As $c_{n+1}$ is defined
 to be zero, all the $c_l$ must vanish. In other words, $E$ is the
 zero element of the cone ${\cal EC}_b(C)$. $\Box$

\medskip

\begin{prop}\label{prop; sim}
 The cone ${\cal EC}_b(C)$ is simplicial. Namely, the generators 
 $e_i, i\in P$ are all linear independent.
\end{prop}

\medskip

\noindent Proof of prop. \ref{prop; sim}
 Suppose that $e_i$ are linearly dependent and there is a linear
 equation $\sum a_i e_i=0$. We can move all the terms with negative $a_i$
 to the right hand side and rewrite the equation as

 $$\sum_{i\in B} b_i e_i=\sum_{j\in B'} b_j e_j, B\cap B'=\emptyset.$$

  This tells us that a single element in ${\cal EC}_b(C)$ has more than
 one expression in terms of the generators $e_i$.

 By prop. \ref{prop; negative}, we can calculate

$$0>(\sum_{i\in B}b_i e_i)\cdot (\sum_{i\in B} b_ie_i)
 =(\sum_{j\in B'}b_j e_j)\cdot (\sum_{i\in B} b_ie_i)\geq 0,$$

 as $B\cap B'=\emptyset$. Contradiction! $\Box$

 The proposition \ref{prop; negative} 
implies that ${\cal EC}_b(C)\cap {\cal EC}_b^{\ast}(C)
=\{ 0 \}$. Namely the dual cone is completely disjoint with the
original exceptional cone.  On the other hand it is easy to see that
 
  $(C-{\cal EC}_b(C))\cap {\cal EC}_b(C)^{\ast}\not=\emptyset$
 if and only if $$({\cal EC}_b(C)-C)\cap (-{\cal EC}_b^{\ast}(C))\not=\emptyset$$

 if and only if $${\cal EC}_b(C)\cap( C-{\cal EC}_b^{\ast}(C))\not= \emptyset $$
 
  We will study the subset of
lattice points in ${\cal EC}_b(C)\cap( C-{\cal EC}_b^{\ast}(C))$
 closely.

\begin{defin}\label{defin; lambda}
 Define the discrete set $\Lambda_b(C)$ to be
 the the lattice points in ${\cal EC}_b(C)\cap( C-{\cal EC}_b^{\ast}(C))$.
\end{defin}

  The elements in $\Lambda_b(C)$ are the expression $\sum_{i\in P}m_ie_i$,
 $m_i\in {\bf Z}^{\geq 0}$,
 such that $(C-\sum_{i\in P}m_ie_i)\cdot e_j\geq 0, j\in P$.

  Let us list the basic properties of the set 
 ${\cal EC}_b(C)\cap( C-{\cal EC}_b^{\ast}(C))$ in the following simple
 proposition,

\begin{prop}
 Let $E(C)$ be the intersection ${\cal EC}_b(C)\cap( C-{\cal EC}_b(C)^{\ast})$
 of two different translated cones. As a proper
 subset of an affine space,
 $E(C)$ is convex as well as unbounded. In fact, given any point $z$ in
 $E(C)$, 
we consider the ray $tz, t\geq 1$, then
 this ray is entirely contained in $E(C)$.
\end{prop} 

\medskip

\noindent Proof: The statement regarding convexity is trivial. 
 Suppose $z\in E(C)$. Then $(C-z)\cdot e_i>0$ for all $i\in P$.
 Then $z\cdot e_i<C\cdot e_i<0$ and we can rewrite 
 $(C-tz)=(C-z)-(t-1)z$ and $(C-tz)\cdot e_i=(C-z)\cdot e_i-(t-1)z\cdot e_i>0$
 for all $i\in P$. Then $tz\in E(C)$. $\Box$

 Because $$2d_{GT}(C)+dim_{\bf R}B=F^2-c_1({\bf K})\cdot F+dim_{\bf R}B+
 \sum_{i\in P} (e_i^2-c_1({\bf K})\cdot e_i)+\Delta_C(E),$$
 
 we are interested at knowing when does $\Delta_C(E)$ take non-positive
 values.

\medskip

\begin{defin}\label{defin; allow}
An lattice element $E$ in $\Lambda_b(C)\subset {\cal EC}_b(C)$ is 
said to be allowable with respect to $C$  
 if the function value $\Delta_C(E)$ is non-positive. 
\end{defin}

 The corresponding decomposition $(C-E, E=\sum m_i e_i)$ is said to be an allowable
 decomposition.

 Let us look at the formula 

$$ \Delta_C(E)=2C\cdot E-E^2-c_1({\bf K})\cdot E-\sum_{i\in P}(e_i^2-
c_1({\bf K})\cdot e_i).$$

 If we view $E$ as a moving variable in $\Lambda_b(C)$ and take an
 element $z$
 in $\Lambda_b(C)$, then $nz,n\geq 1$ form a sequence of lattice points in
 $\Lambda_b(C)$.  Using the fact that the term $E\cdot E$ is negative 
definite, we find that the leading quadratic term in $n$, $n^2z\cdot z$ is always
 positive. Thus it dominates all the other linear or constant 
 terms for the large enough
 $n$. 

 Thus, even though the set $\Lambda_b(C)$ is unbounded, the
 function $\Delta_C(E)$ is bounded below in the un-compact end. As a result,
 the function $\Delta_C(E): \Lambda_b(C)\mapsto {\bf Z}$ 
must attain its absolute minimum somewhere. 

 Before discussing the geometric meaning, let us point out 
 that the locations of the minimums only depend on the numerical
 data and is universally independent to the geometric data of
 the fiber $\pi^{-1}(b)$ itself.

 In principle, the elements $nz$, $n\mapsto \infty$ in the rays will not give
 rise to effective decomposition $C=(C-nz)+nz$ (because
 $(C-nz)\cdot \omega_{{\cal X}/B}\mapsto -\infty$)
 for large enough $n$.  It gives an alternative reason to discard the
 non-compact end of $\Lambda_b(C)$ as one discusses the function
 $\Delta_C(E)$.  However, the non-effectiveness (not being representable by
 pseudo-holomorphic curves) depends on the
 geometric data of the fiber and is not as numerical as the dimension 
 constraint. 

 In principle, the location of the actual minimums relies on the
  data of the various intersection numbers $C\cdot e_i, e_i^2, 
 e_i\cdot e_j$.
 To characterize these
 lattice points geometrically we have to give them a suitable
 interpretation.

 To study the geometric meaning of the minimums of $\Delta_C(E)$, 
let us start with the
 $|P|=1$ case first.  Let us assume that ${\cal EC}_b(C)$ is a one dimensional cone
 generated by a single $e$ with $e^2<0$. 
One suppose that $e$ is represented
 by an irreducible (pseudo)-holomorphic curve above $b\in B$. 
 
 Therefore we write
 $$e^2+c_1({\bf K})\cdot e=2g-2,$$
 where $g$ is the arithmetic genus of $e$ and is usually bigger than the
 geometric genus of $e$ unless the curve dual to $e$ is smooth. 
Any lattice point in ${\cal EC}_b(C)$
 can be written as $m\cdot e, m\in {\bf N}$. Then the function 
 $\Delta_C(E):\{{\bf N} e\}\mapsto {\bf Z}$ can be simplified to

  $$\Delta_C(m)=2m(e\cdot C)-m^2e^2-m(2g-2-e^2)-(2e^2-2g-2).$$

 It is easy to find the value $m_0\in {\bf N}$ 
for which $\Delta_C(m)$ is minimized. We replace
 $m$ by a real variable $x$ and the minimum is
 achieved when the derivative is zero.
 In other words, one has $$2 e\cdot C-2x\cdot e^2-2g+2+e^2=0.$$
 The real number $x$ satisfying this equation is rational. To see the
 geometric meaning of the solution, let us consider the expression
  $e\cdot C-g$, which is always a negative number.

  From simple arithmetic means, it is always possible to re-write
 $eC-g$ as $l\cdot e^2+r$ with $l,r$ non-negative, with $r$ being the remainder
,$0\leq r<-e^2$. It is easy to 
see that $m=l$ is the closest lattice point to the actual minimum of the
 function $\Delta_C(E)$. It is due to the fact that 
 $f(x)=e\cdot C-xe^2-g+1+{e^2\over 2}$ has the property $f(x+1)=f(x)-e^2$ and
 ${e^2\over 2}< f(l)=r+1+{e^2\over 2}\leq {-e^2\over 2}$. 
Therefore $l\cdot e$ is the unique lattice point
 in this range. Moreover it has the crucial property that the function
 $d(m)=\Delta_C(me)$ is monotonically decreasing for $m\geq m_{cri}=l$. 

 In the actual application of the curve counting scheme 
 in the enumerative application, the $g=0$ case is the most interesting situation.

The significant simple property of the function $\Delta_C(E)$ 
in the one dimensional
 case will be used frequently later. Let us summarize it as a lemma.

\begin{lemm}(Moving lemma)\label{lemm; moving}
  Suppose that all the classes $e_i$ satisfy $g(e_i)=0$. 
  Then the function $\Delta_C(E):\Lambda_b(C)\mapsto {\bf Z}$
 is monotonically decreasing if one moves
 from $E$ to $E+e_i$, for all $i$. 
\end{lemm}

\medskip

\noindent Proof: Suppose that $e_i$ has been fixed in our discussion.
First we notice that $E$ can be rewritten as $E=E_0+m_ie_i$ where
 $E_0=\sum_{j\not= i} m_j e_j$.

 Then we compare $\Delta_C(E)$ and $\Delta_{C-E_0}(m_ie_i)$ and see
$$\Delta_C(E)-\Delta_{C-E_0}(m_ie_i)=2C\cdot E-E^2-c_1({\bf K})\cdot E-
  2(C-E_0)\cdot (m_ie_i)+(m_ie_i)^2+c_1({\bf K})\cdot (m_ie_i)$$
$$=2E_0\cdot (m_ie_i)-(E_0+m_ie_i)^2+(m_ie_i)^2-c_1({\bf K})\cdot E_0=
-E_0^2+c_1({\bf K})\cdot E_0,$$ is a constant independent of
 $m_i$.  To show that $\Delta_C(E)$ is monotonically decreasing,
 it suffices to show that $\Delta_{C-E_0}(m_ie_i)$ is monotonically
 decreasing in $m_i$ for $m_ie_i$ in $\Lambda_b(C-E_0)$.  This reduces 
 the problem to the one dimensional case, which has been discussed earlier.
 
 $\Box$

\begin{rem}\label{rem; reallife}
 It should be emphasized that in the single $e$ with $g(e)=0$ case, the integer $l$
 making $\Delta_C(E)$ minimum is the first positive integer $m$ such that 
 $(C-me)$ lies in the dual cone ${\cal EC}_b^{\ast}(C)$. 
 If $C$ is represented by (pseudo) holomorphic curves, then the curve dual to 
 $C$
 can never be irreducible. As $C\cdot e<0$, $C$ must split off as as a
 curve dual to $C-me$ and one dual to $me$ with multiplicity $m$.
 Symbolically we may write $C=(C-me)+me$ as a decomposition of the
 cohomology classes to represent the splitting of curves.
 The only general requirement upon $m$ is that
 $C-me$ and $e$ exist simultaneously as pseudo-holomorphic curves above 
 the fiber of $b$. Thus
 $(C-me)\cdot e\geq 0$ and  it follows that $m\geq {C\cdot
 e \over e^2}\geq l$.  On the other hand, there is no a priori constraint
 about $m$ other than the previous inequality.

 On the other hand, for rational $e$, with $e^2+c_1({\bf K})\cdot e=-2$,
  the minimal choice
 $m=l$ to make $(C-me)\cdot e\geq 0$ 
has the additional nice property that it makes the dimension discrepancy 
function 
$\Delta_C(E)$ minimized. In other words, the curve in the class $C$ tends to split
 off(bubbling off) a certain curves dual to $m$ multiple of $e$ 
such that $C-me$ has nonnegative intersection with $e$. The minimum
 amount $m=l$ also makes the family moduli space dimension 
 $dim_{\bf R}B+d_{GT}(e)+d_{GT}(C-me)$ largest.
 \end{rem}

\medskip

 However this nice topological interpretation does not hold for higher
 genera case. In fact, if the arithmetic genus
  $g$ is larger than $0$, then the role of $e\cdot C$ is replaced by
  $e\cdot C-g$ and it is $e\cdot C-g$  instead of $e\cdot C$ which
 is represented as the form $le^2+r$. Thus the minimum of 
 $\Delta_C(E)$ usually takes place for a larger integer than what 
the naive topological
  constraint predicts.  In other words, the largest expected family dimension  
happens for the multiplicity $m_{cri}$ larger than the topological
 constraint by the integer $[{g\over -e^2}]$ or $[{g\over -e^2}]+1$. \label{genera}

    It is desirable to identify where can the minimum values of 
 $\Delta_C(E)$ occur in $\Lambda_b(C)$.

  Let us consider the translated 
cones $C-{\cal EC}_b(C)^{\ast}+e_i$, the translate of 
 the shifted dual cone $C-{\cal EC}_b^{\ast}(C)$ by the canonical 
basis elements $e_i$. Then we have the
 translated version of the lattice points $\Lambda_b(C)$ in the
 corresponding convex set ${\cal EC}_b(C)\cap C-{\cal EC}_b^{\ast}(C)$, denoted by
 $\Lambda_b(C)_i$. Next we consider the set 

  $$M_b(C)=\Lambda_b(C)-\cup_i(\Lambda_b(C)_i).$$

  The following proposition asserts that $M_b(C)$ is a non-empty finite set.

\begin{prop}\label{prop; finite}
  Let $M(C)$ be the set as defined above, then the set 
$M(C)$ is non-empty and finite.
\end{prop}

\medskip

\noindent Proof:
  Suppose that $M_b(C)$ is empty, then it follows that
 $\Lambda_b(C)\subset \cup_i
 (\Lambda_b(C)_i)$. Let us pick an arbitrary element $\lambda^{(0)}=\lambda$ in the
 left hand side. Then it must be in one of the $\Lambda_b(C)_i$. 
 In other words, $C-(\lambda-e_i)$ also lies in ${\cal C}_E^{\ast}$.
 However this implies that $\lambda^{(1)}=\lambda-e_i$ is also in $\Lambda(C)$. 

 From here we conclude that for any element $\lambda^{(0)}=\lambda$
in $\Lambda(C)$, there must
 be some $i$ such that $\lambda^{(2)}=\lambda-e_i$ still lies in the set $\Lambda_b(C)$. 
 Then by induction, one may use $\lambda-e_i$ instead of $\lambda$ 
 and conclude $\lambda^{(2)}=(\lambda-e_i)-e_j$ is still in $\Lambda_b(C)$.
 However it is
 impossible as it implies that given an element $\lambda^{(0)}=\lambda$, one can
 indefinitely shifts it backward and get another lattice point. Because
 there are only a finite number $e_i$, 
 we must get a lattice point $\lambda^{h}\in \Lambda_b(C)\subset {\cal EC}_b(C)$ 
, $h\gg 0$ 
with at least a negative coordinate entry with respect to some $e_i, i\in P$.
 Contradiction!
 Thus the existence of the lattice point
 in $M_b(C)$ has been derived. 

 To show that $M_b(C)$ is finite, firstly 
we notice that $M_b(C)$ is a discrete subset in
 $H^2(\pi^{-1}(b), {\bf R})$.

\begin{lemm}\label{lemm; cpt}

\end{lemm}
 The set $\overline{{\cal EC}_b(C)\cap (C-{\cal EC}_b^{\ast}(C))\cap 
\bigl(\cup_i(C-{\cal EC}_b^{\ast}(C)
+e_i)\bigr)^c}$ is compact in $H^2(\pi^{-1}(b), {\bf R})$.
\medskip

\noindent Proof of the lemma: 
 Apparently the set is closed, we only need to check it is a bounded 
subset of ${\cal EC}_b(C)\subset H^2(\pi^{-1}(b), {\bf R})$.

  Suppose that $E\in \overline{{\cal EC}_b(C)\cap (C-{\cal EC}_b^{\ast}(C))\cap 
\bigl(\cup_i(C-{\cal EC}_b^{\ast}(C)
+e_i)\bigr)^c}$, then $E\in {\cal EC}_b(C)\cap (C-{\cal EC}_b^{\ast}(C))$
 but $E$ is not in the interior of 
${\cal EC}_b(C)\cap (C-{\cal EC}_b^{\ast}(C)+e_i)$ for 
 all $i\in P$.

  Thus, $(C-E)\cdot e_i\geq 0, i\in P$ but $(C-E+e_j)\cdot e_i\leq 0, i, j\in P$.

 This implies that the values $E\cdot e_i, i\in P$ are bounded by

 $$min_{j\in P}(C+e_j)\cdot e_i\leq  E\cdot e_i\leq C\cdot e_i, i\in P.$$

  This implies that the pairing functionals $e_i\cdot \circ$ take bounded
 values on $E$. On the other hand, $e_i, i\in P$ form a basis of the
 minimal vector space containing ${\cal EC}_b(C)$, the pairing functional
 $e_i\cdot \circ$ is a linear coordinate system. As $E$ has bounded 
coordinates, such $E$ forms a bounded subset in ${\cal EC}_b(C)$. $\Box$

 From the fact that
 $\overline{{\cal EC}_b(C)\cap (C-{\cal EC}_b^{\ast}(C))\cap 
(\cup_i(C-{\cal EC}_b^{\ast}(C)
+e_i))^c}$ is compact and $M_b(C)$ is a discrete subset,
 it follows that $M_b(C)$ has to be a finite set. $\Box$

 Usually it is not clear whether the
 lattice points in the finite set $M_b(C)$ are unique or not. 
 The following proposition clarifies the significance of the lattice points
 in $M_b(C)$.

\begin{prop}\label{prop; max}
 Let $E=\sum_i m_i e_i, m_i\in {\bf N}$ be in 
$\Lambda_b(C)\subset {\cal EC}_b(C)\cap (C-{\cal EC}_b^{\ast}(C))$.
 Then there exists at least one element $z$ in $M_b(C)$ 
 and a finite sequence of elements $z_p\in \Lambda_b(C)$, $1\leq p\leq N$
 such that

\medskip

(i). $z_1=z$, $z_N=E$.

(ii). If $p\not=1$, then $z_p-z_{p-1}=e_i$ for some $1\leq i\leq |P|$.

\end{prop}

\medskip

\noindent Proof of the proposition: $\Box$
 Given an element $E\in \Lambda_b(C)$, if $E$ also lies in $M_b(C)$, we
 are done. We set $N=1$ and $z_1=E=z\in M_b(C)$.

 Otherwise if $x_0=E\not\in M_b(C)$, then $E\in \Lambda_b(C)_i$ 
for some $1\leq i\leq |P|$. In other words, there exists an $x_1\in \Lambda_b(C)$
 such that $x_0=E=x_1+e_i$ for some $1\leq i\leq |P|$.

 Take $x_1\in \Lambda_b(C)$ and repeat the argument, either one gets 
 $x_2\in \Lambda_b(C)$ such that $x_1=x_2+e_i$ for some $1\leq i\leq |P|$
 or $x_2\in M_b(C)$. By induction, one may get a sequence of points 
 $x_n\in \Lambda_b(C)$. One argues that for a large enough $x_n$,
  we must have $x_n\in M_b(C)$, or the induction process never stops and
 one gets an infinite sequence $x_n\in \Lambda_b(C), n\in {\bf N}$. Because
 there are a finite number of $e_i, 1\leq i\leq |P|$, 
such an $x_n=\sum_{1\leq i\leq |P|}\underline{m}_ie_i$ must has a negative
 entry $m_j<0$ for some $j\leq |P|$ and
 falls out of ${\cal EC}_b(C)$. Contradiction!

 Then we take $z_1=z=x_n\in M_b(C)$ and rename the sequence $x_p, 1\leq p\leq n$
 by $z_q=x_{n+1-q}, q\leq n+1=N$.

 Such a finite sequence of lattice elements satisfies 

 (i). $z_1=z\in M_b(C)$, $z_N=x_0=E\in \Lambda_b(C)$. $z_q\in \Lambda_b(C)$

 (ii). For $p>1$, $z_p-z_{p-1}=e_i$ for some $i\leq |P|$.

\medskip
 
 We are done. $\Box$

\medskip

 We assume that $g(e_i)=0$ for all $1\leq i\leq |P|$ in the following
 remark.

\medskip

\begin{rem} \label{rem; geometry}
  Suppose that $C$ is represented by a pseudo-holomorphic curve over the
 point $b\in B$, then one
 argues that the curve must contain irreducible curves dual to the various
 $e_i, i\leq |P|$. It is because $e_i$ is known to be irreducible 
and pseudo-holomorphic above $b\in B$, then $C\cdot e_i\geq 0$ if all the
 irreducible components of the curve dual to $C$ are distinct from the one
 dual to $e_i$.

 Thus, one may write 
$C=(C-\sum_i m_i e_i)+(\sum_i m_i e_i)$ with $m_i\in {\bf N}$, where
 $(C-\sum_i m_i e_i)$ is dual to a curve disjoint from all the exceptional 
 curves dual to $e_i$. Then we must have $(C-\sum_i m_i e_i)\cdot e_j\geq 0$
 for all $j\leq |P|$. If we take $E=\sum_i m_i E_i\in {\cal EC}_b(C)$,
 we must have $(C-E)\in {\cal EC}_b^{\ast}(C)$. In other words,
 $E\in \Lambda_b(C)\subset {\cal EC}_b(C)\cap \bigl(C-{\cal EC}_b^{\ast}(C)\bigr)$.

We have use an $E$ to resemble the pattern $C$ splits into different multiples
 of $e_i, i\leq |P|$.

 If one imagines a pseudo-holomorphic curve dual to $C-E_0$ splits into a curve
 dual to $C-E_0-e_i$ and a curve dual to $e_i$ as a degeneration process (known
 as the bubbling off phenomenon in symplectic geometry), then
 the lattice move $E_0\mapsto E_0+e_i$ in $\Lambda_b(C)$ is equivalent to
 ``bubbling off'' an unit of $e_i$ from $C$.

 The path from $z$ to $E$ by a finite sequence of lattice moves as in the
 condition $(ii)$ of prop. \ref{prop; max} indicates that the curve dual to
 the sum of $(C-E)$ and $E$ can be viewed formally as degenerated from a curve
 dual to $C-z$ and $z$ by a finite number of bubbling offs of $e_i, i\leq |P|$.

\medskip

(a). In other words, if we consider the moduli space (and its bubbling off) 
of curves dual to $C$ splitting
 into one dual to $C-z$ and a combination of multiple coverings of exceptional 
 curves dual to $z$, it contains the given curve dual to the sum of $C-E$ and $E$.

\medskip

(b). Given any pseudo-holomorphic curve dual to $C$, one may assign an unique
 $E\in \Lambda_b(C)$ to it. By prop. \ref{prop; max}, one may associate
 the point $z\in M_b(C)$ to $E$ with the properties,

\medskip

 (i). The decomposition $C$ into $(C-z)+z$ has the highest expected
 family Gromov-Taubes dimension among all
 the decompositions into distinct pseudo-holomorphic curves
 related by the elementary moves.

\medskip

 (ii). The curve dual to $(C-E)+E$ can be thought as a degenerated version of
 curves dual to $(C-z)+z$ by bubbling off a few rational exceptional curves
 dual to the $e_i$.
\end{rem}

 When one applies the family switching formula of rational exceptional curves,
 one is able to rewrite some mixed family invariant of $C-\sum e_i$ over
 a locus over which all the exceptional curves in $e_i, i\leq |P|$ co-exist, (
 with all the multiplicities $m_i\equiv 1$) in terms of the family invariants
 of $C-z$, $z=\sum_{1\leq i\leq |P|} n_i e_i\in M_b(C)$. Combine with remark
 \ref{rem; geometry}, this tells us that formally the mixed family
 invariant can be interpreted as the counting of curves in the class
 $C-\sum_{1\leq i\leq |P|}n_ie_i$.

\medskip

 Returning to the pointwise discussion over $b\in B$, 
any two different lattice points within $M_b(C)$ cannot be related
 by each other by effective translation. i.e.\ shifting from one to another by
 $\sum_i c_i e_i,c_i\leq 0$.  In reality the geometric meaning of
 shifting toward the right means degeneration or bubbling off phenomena(in the
 case $g(e_i)=0$). In this way, we may give the lattice points in
 $\Lambda_b(C)$ a partial ordering and different elements in $M_b(C)$ are the
 maximum elements (not necessarily the greatest element) of the partial ordering. 

\begin{defin}
 Let $\lambda_1$ and $\lambda_2$ be two different
  lattice elements in $\Lambda_b(C)$. 
We say that $\lambda_1$ is greater than $\lambda_2$, denoted by $\lambda_1 
 \sqsupset \lambda_2$
 if $lambda_2$ can be gotten from $\lambda_1$ by a finite sequence
 of effective lattice translations, i.e.\ there exists
 $\lambda_2=z_n, z_{n-1}, z_{n-2}, \cdots, \lambda_1=z_1 \in \Lambda_b(C)$ such 
that $z_p-z_{p-1}=n_pe_{i_p}$ for some $n_p\in {\bf N}$, $1\leq i_p\leq |P|$.
  It is apparent that 
 this relation $\sqsupset$ is transitive. 

We say that an element $\lambda$ is a maximal
 element if there is no other element in $\Lambda_b(C)$ which is greater than
 it under the partial ordering $\sqsupset$. 
\end{defin}

\medskip

  Then the set $M_b(C)$ consists of the maximum elements in
 $\Lambda_b(C)$. This justifies the notation $M_b(C)$ and its dependence on
 the class $C$.  

 From the previous lemma \ref{lemm; moving}
 one finds that if $\lambda_1\sqsupset \lambda_2$, then
 $\Delta_C(\lambda_1)<\Delta_C(\lambda_2)$. Therefore, in the rational case,
 $M_b(C)$ also represents the lattices points of $\Lambda_b(C)$
 whose $\Delta_C(E)$ values
 are smallest. Even though the values of $\Delta_C(E)$ are different for
 the  different
 lattice points, it is no use to compare their value as they are not linked
 to each other by effective lattice shifting. On the other hand, our
 discussion does not rule out the possibility that by different
 effective shiftings one can reach from two different maximum elements $\in
 M_b(C)$ to the same lattice point $\lambda\in \Lambda_b(C)$.
 As effective shiftings correspond to pseudo-holomorphic curve
 degenerations, the single decomposition $(C-\lambda, \lambda)$, 
$\lambda\in \Lambda_b(C)$ can possibly be
 degenerated from two distinct maximal decompositions 
 $(C-\lambda_1, \lambda_1), (C-\lambda_2, \lambda_2)$, $\lambda_1, \lambda_2
 \in M_b(C)$, $\lambda_1\sqsupset \lambda, \lambda_2\sqsupset \lambda$. 

 In the latter part of the paper,
 we would like to discuss how does the degeneration process affect the 
relative
 obstruction bundles. 
 The family switching formula assigns a relative obstruction bundle
 to a degeneration of  decompositions such that 
$\Delta_C(E)-\Delta_C(E+\sum n_ie_i)$ is directly related to the rank of
the bundle.
It turns out that not only on the numerical level do 
 the partial ordering organizes the  lattice points
 of $\Lambda_b(C)$ in a nice way, but they also co-relate
  the relative obstruction bundles. 
 This will be the main focus of the next section.

\medskip

\subsection{\bf The Irrational, $g(e_i)>0$, Cases} \label{subsection; irr}

\medskip

 In the previous discussion, we have focused upon the $g(e_i)=0$ case. 
Let us consider the general situation that 
arithmetic genera $g(e_i), i\leq |P|$
 may be non-zero.  As we discuss earlier (see page \pageref{genera}), 
the appearance of the genus term 
shifts the minimum value of $\Delta_C(E)$. This happens
 even when the exceptional effective cone is of one dimension.

 Let us make this simplifying assumption temporally.
 Once we have a closer look at the formula
 $-x\cdot e^2+e\cdot C-g+1+e^2/2$, 
 we may collect $e\dot C-g+1$ into a single term with a particular topological
 meaning. Imagining that $e$ is represented by a smooth pseudo-holomorphic curve in
 $\pi^{-1}(b)$. Suppose $C$ is the first Chern class of a holomorphic line
 bundle ${\bf F}$ on the Riemann surface $\Sigma\subset \pi^{-1}(b)$.
 Then the expression $e\cdot C-g+1=\int_{\Sigma}c_1({\bf F})-g(\Sigma)+1$
 resembles the holomorphic Euler number of the holomorphic line bundle ${\bf F}$.
 The appearance of the special number indicates that the tangent-obstruction
 complex would contain a term associated to 
$H^0(\Sigma, {\bf F})- H^1(\Sigma, {\bf F})$, 
 $e=PD[\Sigma]$.

 In general, consider the intersection matrix $I_{i,j}=e_i\cdot
 e_j$, $1\leq i, j\leq |P|$. 
 Let us take the dual basis $e_i^{\ast}$ with respect to $e_i$ such that
 $e_i\cdot e^{\ast}_j=I_{i,j}\delta_{i,j}$,
 and the normalized dual basis ${e^{\ast}_j\over
 e_j\cdot e_j}=\hat {e}_j$. Then we have $e_i\cdot \hat {e}_j=\delta_{i,j}$ for
 $1\leq i, j\leq |P|$. 
 
 As the classes $C$ and $c_1({\bf K}_{{\cal X}/B})$ 
may not lie in the cone ${\cal EC}_b(C)$, we 
 project them orthogonally into the subspace  ${\cal EC}_b(C)
\otimes {\bf R}$. 
 Denote ${\hat C}:=\sum_{1\leq i\leq |P|}(C\cdot e_i){\hat e_i}$ and 
  ${\hat K}:=\sum_{1\leq i \leq |P|}(c_1({\bf K}_{{\cal X}/B})\cdot e_i)
{\hat e_i}$, then we must have

$$c_1({\bf K}_{{\cal X}/B})\cdot e_j={\hat K}\cdot e_j, 
C\cdot e_j={\hat C}\cdot e_j.$$
The ${\hat C}$ and ${\hat K}$ lie in the subspace ${\cal
   EC}_b(C)\otimes {\bf R}
\subset H^2(M, {\bf R})$ and are the projection of $C$ and 
 $c_1({\bf K}_{{\cal X}/B})$ to the subspace.  These elements
 ${\hat C}$ and ${\hat K}$ are rational points of the minimal subspace 
 ${\cal EC}_b(C)\otimes {\bf R}$. 
   
 Differentiating $\Delta_C(E)|_{E=\sum m_i e_i}$ 
with respect to the variable $m_i$, one derives that 

 $$-2E\cdot e_i+2{\hat C}\cdot e_i-{\hat K}\cdot e_i=0,1\leq i\leq |P|.$$

 Then one concludes that the rational point $E=-{\hat C}+{1\over 2}{\hat K}$
 
 is the minimum of the quadratic function $\Delta_C(E)$. 
  By using the adjunction equality one may rewrite
 $${\hat K}=\sum_i(2g(e_i)-2-e_i^2){\hat e_i},$$
 $${\hat C}-{1\over 2}{\hat K}=\sum_i(e_i\cdot C-g(e_i)
+1+{e_i^2\over 2}){\hat e_i}.$$ 

  If the arithmetic genera $g(e_i)$ 
change from zero to positive values, then we find the
 minimum of the function $\Delta_C(E)$ are shifted from their original value by
 the rational vector 
 $-g(E):=-\sum_{1\leq i \leq |P|}g(e_i){\hat e_i}$. It is easy to see that
 $g(E)\in {\cal EC}_b(C)\otimes {\bf R}$
 is a rational element in ${\cal C}_E^{\ast}$;
  the dual cone ${\cal EC}_b^{\ast}(C)$. As we have proved (in prop. 
\ref{prop; negative}) that
 ${\cal EC}_b(C)\cap {\cal EC}_b^{\ast}(C)=\{0 \}$,
 so the element $g(E)$ never lies in 
the original ${\cal EC}_b(C)$.  
 
 But it is not clear from the definition whether $-g(E)$ can be in
 ${\cal EC}_b(C)$ or not.

\medskip

\section{\bf The Global Discussion and The Admissible Decomposition Classes}
\label{section; decomposition}

\bigskip

 In the previous section, we have finished the pointwise discussion on the
 exceptional cones of $C$ and the family dimensions.
  We would like to patch the local discussion over $b\in B$ 
together and introduce some additional structure on the base space $B$ of
 the fibration $\pi:{\cal X}\mapsto B$.
 
\medskip

Suppose that $e$ is an exceptional class in the sense of definition 
\ref{defin; exception}. 
Suppose that at a given point $b\in B$ the class $e$ has been
 represented by an irreducible holomorphic curve in ${\cal X}_b$, one may
consider the locus $S_e\subset B$ over which the class $e$ is effective.
 It is well known that $S_e$ must be a compact subset of the compact set $B$.

 By adjunction equality we have 
$e^2+c_1({\bf K}_{{\cal X}/B})\cdot e=2g_{arith}(e)-2$. On the other hand, the
 Gromov theory predicts that the 'expected' dimension of the set $S_e$ is
 $dim_{\bf C}B+{e^2-c_1({\bf K}_{{\cal X}/B})\cdot e\over 2}$.

In order for the curve to be generic, i.e. expected dimension of 
$S_e\geq dim_{\bf C}B$, both inequalities 
$$e^2-c_1({\bf K}_{{\cal X}/B})\cdot e\geq 0, 
e^2+c_1({\bf K}_{{\cal X}/B})\cdot e\geq -2$$

  have to be satisfied. In particular, $e^2=-1$ and $e$ is represented by
 a genus zero smooth curve. Thus, we may conclude that the only generic 
 exceptional curve within a family $\pi:{\cal X}\mapsto B$ has to be
 a smooth $-1$ curve.

 For simplicity, let us assume that the arithmetic genus of the curve dual 
 to $e$, 
$g_{arith}(e)=0$, in the following
 discussion. Suppose that $e^2=-n$, 
 the irreducible rational curve representing $e$ is called a $-n$ rational 
(exceptional) curve. According to the dimension formula, the expected 
complex dimension of $S_e$ is $dim_{\bf C}B+{d_{GT}(e)\over 2}=dim_{\bf C}B+1-n$.
 This indicates when the self-intersection number $e\cdot e$ is more negative,
 the expected dimension of $S_e$ would be much lower.

\medskip

  In general an irreducible holomorphic curve may degenerate and breaks into
 more than one irreducible component. Then there is a subset 
$S_e^{sm}\subset S_e$ over which the curve representing $e$ is 
smooth and irreducible. The set $S_e-S_e^{sm}$ is the locus over which
 the curve representing $e$ breaks into more than one component. In an
 ideal situation, $S_e^{sm}$ is dense in $S_e$. Then $S_e-S_e^{sm}$ can be
 thought as the ``boundary points'' of $S_e^{sm}$ collecting all the
 degenerated configurations of $e$.

\medskip

 When more than one $e$, say $e_1, e_2, \cdots, e_k$ are represented by
 holomorphic curves at the same $b\in B$, the set of all such $b$ is 
nothing but $S_{e_1}\cap S_{e_2}\cdots\cap S_{e_k}=\cap_{i\leq k}S_{e_i}$.

 In terms of intersection theory, the expected dimension of $\cap_{i\leq k}S_{e_i}$
 is $dim_{\bf C}B+{1\over 2}\sum_{i\leq k}d_{GT}(e_i)$.

\medskip

 Let us list all the basic assumptions on a ``perfect'' 
set of $C$-exceptional classes.

\medskip

\begin{defin}\label{defin; perfect}
 Let $C$ be a $(1, 1)$ class on ${\cal X}$ which restricts to non-trivial
 class on an algebraic fibration $\pi:{\cal X}\mapsto B$.

 A finite set of exceptional classes $Q$ on $\pi:{\cal X}\mapsto B$ is said to 
be perfect if they satisfy the following list of basic assumptions.

\begin{assum}\label{assum; perfect}

\noindent (i). $e\cdot C<0$.

\noindent (ii). Either $S_e=\emptyset$ or $S_e$ is a
 $dim_{\bf C}B+{1\over 2}d_{GT}(e)$ dimensional closed sub-variety of $B$.

\noindent (iii). The set $S_e^{sm}$ is dense in $S_e$ consisting smooth points
 of $S_e$.

\noindent (iv). If ${\bf e}$ is a holomorphic curve in ${\cal X}_b, b\in B$
 representing $e$, then all the irreducible components of ${\bf e}$ represent
 exceptional classes. That is to say,
 none of the irreducible components has a non-negative
 self-intersection number. Moreover, if this exceptional class also
 satisfies (i). then
 it must also be in the original collection $Q$.

\noindent (v). Let ${\bf e}$ be a smooth irreducible holomorphic curve
 in ${\cal X}_b$ representing $e$ and $b\in S_e^{sm}$. Then the 
 composite morphism 

$${\bf T}_bB\mapsto H^1({\cal X}_b, \Theta{\cal X}_b)
\mapsto H^1({\cal X}_b, {\cal O}_{{\cal X}_b}({\bf e}))\mapsto H^1({\bf e}, 
{\cal N}_{\bf e}{\cal X}_b)$$

 induces an isomorphism 
${\bf N}_{S_e^{sm}}B|_b\mapsto H^1({\bf e}, {\bf N}_{\bf e}{\cal X}_b)$.

\noindent (vi). Suppose $e_1, e_2, \cdots, e_k$, $e_i\cdot e_j\geq 0$, $i\not= j$
 are within the
collection $Q$ and are satisfying all (i).-(v)., then the locus of co-existence of
$e_1, e_2, \cdots, e_k$, $\cap_{i\leq k}S_{e_i}$, is either empty or
 is a $dim_{\bf C}+{1\over 2}\sum_{i\leq k}d_{GT}(e_i)$ sub-variety in $B$.

\noindent (vii). The set $\cap_{i\leq k}S_{e_i}^{em}$ is dense
 in $\cap_{i\leq k}S_{e_i}$ containing the smooth points in it.
\end{assum}
\end{defin}

\medskip

  Under the condition (ii)-(iii) in the assumption 
\ref{assum; perfect}, each $S_e$ defines an algebraic cycle class $[S_e]\in 
{\cal A}_{dim_{\bf C}B+{1\over 2}d_{GT}(e)}(B)$ and can be
 used to define the algebraic family Seiberg-Witten invariant of $e$.
 The cycle class $[S_e]$ can be viewed as the moduli cycle of the 
 exceptional class $e$.

  Likewise, the locus of co-existence $\cap_{i\leq k}S_{e_i}$ of various
 $e_i, 1\leq i\leq k$, also
 defines an algebraic cycle class $[\cap_{i\leq k}S_{e_i}]$, which is 
identical to the intersection cycle class $\cap_{i\leq k}[S_{e_i}]\in 
{\cal A}_{dim_{\bf C}B+{1\over 2}d_{GT}(e_i)}(B)$. 

\medskip

\begin{rem}\label{rem; v}
 In condition (v), the sheaf cohomology $H^1({\cal X}_b, \Theta{\cal X}_b)$
 parametrizes the infinitesimal Kodaira-Spencer 
deformations of complex structures on ${\cal X}_b$. Then
 ${\bf T}_bB\mapsto H^1({\cal X}_b, \Theta{\cal X}_b)$ is the tautological
 map induced by the infinitesimal deformation of the family of algebraic surfaces 
${\cal X}\mapsto B$.
 The morphism $H^1({\cal X}_b, \Theta {\cal X}_b)\mapsto 
H^1({\bf e}, {\cal N}_{\bf e}{\cal X}_b)$ is induced by the splitting
 $\Theta {\cal X}_b|_{\bf e}={\cal T}_{\bf e}\oplus {\cal N}_{\bf e}{\cal X}_b$
along the holomorphic curve ${\bf e}\subset {\cal X}_b$.
\end{rem}

\medskip

Let $b\in B$ and let ${\cal EC}_b(C)$ be as defined in definition \ref{defin; exc}.
By proposition \ref{prop; sim} the cone ${\cal EC}_b(C)$ is a simplicial cone
 generated by a collection of $e_i$ represented by irreducible
 exceptional curves in ${\cal X}_b$.

\begin{defin}\label{defin; Qcone}
 Define ${\cal EC}_b(C; Q)$ be the simplicial sub-cone of 
 ${\cal EC}_b(C)$ generated by the elements $e_i\in Q$.
\end{defin}

\medskip

 As $b$ moves on $B$, the cone ${\cal EC}_b(C; Q)$ often 
changes along with $b$.
 It makes sense to ask the following question,

\medskip

\noindent {\bf Question:} Describe the pattern of the variations of
 ${\cal EC}_b(C)$, $b\in B$, in terms of algebraic geometric datum on $B$.

\medskip

Let ${\cal C}$ be a simplicial cone generated by elements in $Q$.

\begin{defin}\label{defin; stratum}
Define $S_{\cal C}\subset B$ to be the set of all $b\in B$ such that
 ${\cal EC}_b(C; Q)\equiv {\cal C}$. 
\end{defin}

For certain ${\cal C}$ there can be no such $b\in B$ and $S_{\cal C}$ is empty.
We are only interested at those ${\cal C}$ with a non-trivial $S_{\cal C}$ and
 consider the pair $(S_{\cal C}, {\cal C})$.
\medskip

\begin{prop}\label{prop; const}
Let $Q$ be a perfect finite set of exceptional classes. 
Let $(S_{\cal C}, {\cal C})$ be a pair with $S_{\cal C}\not=\emptyset$.
 Then $S_{\cal C}$ is a locally closed subset of $B$ and 
$S_{\cal C}\subset \overline{S_{\cal C}}$ consists of smooth points of
 $\overline{S_{\cal C}}$.
\end{prop}

\noindent Proof of proposition \ref{prop; const}:
 Suppose that ${\cal C}$ is generated by $e_i\in Q, 1\leq i\leq k$.
 Because $Q$ is perfect, by assumption \ref{assum; perfect} (ii)., (iii)., 
(vi)., (vii)., $\cap_{1\leq i\leq k}S_{e_i}$ is a 
 $dim_{\bf C}B+{1\over 2}\sum_{i\leq k}d_{GT}(e_i)$ dimensional sub-variety in $B$.

 By our definition of ${\cal EC}_b(C; Q)$ as a subcone of ${\cal EC}_b(C)$, each 
$e_i\in {\cal C}\equiv {\cal EC}_b(C; Q)$ has to be represented by a smooth and
 irreducible exceptional curve above $b$. Thus, 
$S_{\cal C}\subset \cap_{1\leq i\leq k}S_{e_i}^{sm}$.
 The smoothness of $S_{\cal C}$ follows from condition (vii) of
 assumption \ref{assum; perfect}.

 We plan to argue that $S_{\cal C}$ is open and dense
 in the closed set $\cap_{1\leq i\leq k}S_{e_i}$. Then
 $\overline{S_{\cal C}}=\cap_{1\leq i\leq k}S_{e_i}$.
 Once this is achieved, the local closeness of $S_{\cal C}$ follows.

 Firstly, by deformation theory of smooth curves, the set $S_{e_i}^{sm}$ is
 open in $S_{e_i}$. Thus $\cap_{1\leq i\leq k}S_{e_i}^{sm}$ is open in
 $\cap_{1\leq i\leq k}S_{e_i}$.

 Denote the difference $\cap_{1\leq i\leq k}S_{e_i}^{sm}-S_{\cal C}$ as $A$.
 We argue that its closure $\overline{A}$ is a higher codimension subset in
 $\cap_{1\leq i\leq k}S_{e_i}$.

 At all $b\in A$, the class $e_i$, $1\leq i\leq k$ are
 represented by irreducible exceptional curves. On the other hand,
 $A\cap S_{\cal C}=\emptyset$. Thus, the cones ${\cal EC}_b(C; Q)$, $b\in A$
 have to be strictly larger than ${\cal C}$. In other words, for every
 $b\in A$, there must
 be some additional exceptional classes $\tilde{e}\in Q$, effective and
 irreducible over $b$. One may collect all such $\tilde{e}\in Q\in
 {\cal EC}_b(C; Q)-{\cal C}$ into a set $E$ when $b$ runs through
 all the points in  $A$.

 It is apparent that $A=\cap_{1\leq i\leq k}S_{e_i}^{sm}-S_{\cal C}\subset
 \cup_{\tilde{e}\in E}(\cap_{1\leq i\leq k}S_{e_i}\cap S_{\tilde{e}})$.

 By condition (vii) of assumption \ref{assum; perfect}, 
 $(\cap_{1\leq i\leq k}S_{e_i}\cap S_{\tilde{e}})$ is of dimension
 $dim_{\bf C}B+{1\over 2}(\sum_{1\leq i\leq k}d_{GT}(e_i)+d_{GT}(\tilde{e}))$ and
 is of lower dimension than $\cap_{1\leq i\leq k}S_{e_i}$ because 
 $d_{GT}(\tilde{e})<0$.

Because the set $E$ is finite,
 $\cup_{\tilde{e}\in E}(\cap_{1\leq i\leq k}S_{e_i}\cap S_{\tilde{e}})$,
 a finite union of closed subsets of $\cap_{1\leq i\leq k}S_{e_i}$ is
 closed.

Thus, one may write $S_{\cal C}$ as 

$$\cap_{1\leq i\leq k}S_{e_i}^{sm}-
\cup_{\tilde{e}\in E}(\cap_{1\leq i\leq k}S_{e_i}\cap S_{\tilde{e}})$$

$$=\cap_{1\leq i\leq k}S_{e_i}-(\cap_{1\leq i\leq k}S_{e_i}-
\cap_{1\leq i\leq k}S_{e_i}^{sm})-
\cup_{\tilde{e}\in E}(\cap_{1\leq i\leq k}S_{e_i}\cap S_{\tilde{e}}),$$

and apparently is
 open in $\cap_{1\leq i\leq k}S_{e_i}$. $\Box$

\medskip

 Because $Q$ is a finite set, there are finitely many possible ${\cal C}$ with
 non-empty $S_{\cal C}$. By its definition, 
$S_{{\cal C}_1}\cap S_{{\cal C}_2}=\emptyset$ if ${\cal C}_1\not={\cal C}_2$.

\medskip

\begin{prop}\label{prop; closure}
Let $(S_{\cal C}, {\cal C})$ be a pair with $S_{\cal C}\not=\emptyset$.
The boundary set of $S_{\cal C}$, $\overline{S_{\cal C}}-S_{\cal C}$,  
is contained inside
 a disjoint union of different $S_{{\cal C}'}$ with ${\cal C}\subset {\cal C}'$.
\end{prop}

\medskip

\noindent Proof of proposition \ref{prop; closure}: Take an arbitrary
 $b\in \overline{S_{\cal C}}-S_{\cal C}$, there are three different
 possibilities.

\medskip

\noindent (A). All the $e_i$, $1\leq i\leq k$ which generates ${\cal C}$
 still represent irreducible curves in ${\cal X}_b$. But some other new
 exceptional class in $Q$ becomes effective over $b$ and becomes
 a generator of ${\cal EC}_b(C; Q)$.

\medskip

\noindent (B). The curves representing 
some or all of the $e_i$, $1\leq i\leq k$ break into
 more than one irreducible components.

\medskip

\noindent (C). Some new exceptional class $\tilde{e}$, $\tilde{e}\cdot e_i\geq 0$
 becomes effective over $b$, while the curves representing some of the $e_i$
 break into more than one component. (The mixture of (A). and (B).)

\medskip

 Suppose that $\tilde{e}_j\in Q, 1\leq j\leq l$ are the new exceptional class(es) 
in case (A)., consider the cone ${\cal C}'$ generated by $e_i, 1\leq i\leq k$ and
 $\tilde{e}_j, 1\leq j\leq \tilde{k}$. 
Such points $b$ are contained in the set $S_{{\cal C}'}$ and apparently
 ${\cal C}\subset {\cal C}'$.

 In case (B)., suppose that, say $e_1$, has broken into components,
$e_1=\sum_{1\leq r\leq r_1} e_{1;r}$. Then for at least one $e_{1; r}$,
 $e_{1; r}$ pairs negatively with $C$, i.e. $e_{1; r}\cdot C<0$. If not, 

$$0>e_1\cdot C=(\sum_{1\leq r\leq r_1}e_{1;r})\cdot C=\sum_{1\leq r\leq r_1}
e_{1; r}\cdot C\geq 0,$$ 
contradicting to the condition (i). of the perfectness assumption.
On the other hand, the condition (iv) of perfectness assumption 
\ref{assum; perfect} implies that $e_{1;r}\in Q$ if $e_{1;r}\cdot C<0$.
 So ${\cal EC}_b(C; Q)\not=\emptyset$.

 Such a $b$ is in $S_{{\cal C}'}$ with ${\cal C}'={\cal EC}_b(C; Q)$. 
Because $b$ is in the closure of
 $S_{\cal C}$ over which the classes in ${\cal C}$ are 
effective, by the degeneration argument 
all classes in ${\cal C}$ remain effective over $b$ as well.
 This implies that ${\cal C}\subset {\cal EC}_b(C; Q)={\cal C}'$.

The discussion for the case (C). is similar and we leave it to the reader. $\Box$

\medskip

 The proposition motivates us to define a partial ordering among
 different $(S_{\cal C}, {\cal C})$.

\medskip

\begin{defin}\label{defin; partial}
 The pair $(S_{\cal C}, {\cal C})$ is said to be greater than
 $(S_{\cal C}',  {\cal C}')$ under $\succ$, denoted as
$(S_{\cal C}, {\cal C})\succ (S_{\cal C}',  {\cal C}')$,
 if  ${\cal C}\subset {\cal C}'$.
\end{defin}

Notice that $\succ$ is a necessary condition for $\overline{S_{\cal C}}\cap
 S_{{\cal C}'}$ to be non-empty.

\medskip

\subsection{\bf The Admissible Decomposition Classes over $S_{\cal C}$}
\label{subsection; adm}

\bigskip

 Having addressed the structure of $S_{\cal C}$, we move ahead to
 address the decomposition classes.

 Consider the disjoint union $\coprod_{{\cal C}, S_{\cal C}\not=\emptyset}
 S_{\cal C}$. Either it is equal to the whole $B$, or it is a closed 
subset (by proposition \ref{prop; closure}) of $B$.
 In the first case, some $S_{\cal C}$ is of top dimension in $B$ and
 ${\cal C}$ is generated by a finite number of $-1$ classes in ${\cal C}$.

 If it is the case, any effective curve dual to $C$ in the fibers of the
 family ${\cal X}\mapsto B$ over $S_{\cal C}$
  must break off a certain multiples of $-1$ curves and the family
 theory suffers the same symptom as in Mcduff's proposal [Mc] ($B=pt$).

 The more interesting situation is when 
$\coprod_{{\cal C}, S_{\cal C}\not=\emptyset} S_{\cal C}\subset B$
 is a proper closed subset.

 In such a situation, any effective curve dual to $C$ over
 $b\in \coprod_{{\cal C}, S_{\cal C}\not=\emptyset} S_{\cal C}$ has to
  break off certain multiples of exceptional curves in ${\cal EC}_b(C; Q)$.
 What types of exceptional curves it has to break off depends on which
 $S_{\cal C}$ does $b$ lie in. 

 \medskip

 A sketch of the general enumerative application of our curve counting scheme
 to the family algebraic Seiberg-Witten invariants 
 ${\cal AFSW}_{{\cal X}\mapsto B}(1, C)$ can be outlined as the following.

 (1). Break the whole base space $B$ into different $S_{\cal C}$ and
 $B-\coprod_{{\cal C}; S_{\cal C}\not=\emptyset} S_{\cal C}$.

 This set level decomposition gives a stratification of $B$ and 
$B-\coprod_{{\cal C}, S_{\cal C}\not=\emptyset} S_{\cal C}$ is the only
 top dimensional stratum. The dimension of each $S_{\cal C}$ can 
 be calculated through the argument of proposition \ref{prop; const}.

\medskip

 (2). Attach a local family invariant contribution to each 
$\overline{S_{\cal C}}
\not=\emptyset$ depending on the generators $e_i, 1\leq i\leq k$
 of ${\cal C}$ and the breaking of $C$ into $C-\sum_{1\leq i\leq k}e_i$
 and $\sum_{1\leq i\leq k}e_i$. (see [Liu1] and [Liu3] for
 some examples)
 The invariant contribution is some mixed family invariant of the form
${\cal AFSW}_{{\cal X}\times_B \overline{S_{\cal C}}\mapsto \overline{S_{\cal C}}}
(\cdot, C-\sum_{1\leq i\leq k}e_i)$.

 We expect the local family invariant contribution to be nonzero only if 
 $dim_{\bf C}B+{1\over 2}
d_{GT}(C-\sum_{1\leq i\leq k}e_i)+{1\over 2}\sum_{1\leq i\leq k}d_{GT}(e_i)\geq
dim_{\bf C}B+d_{GT}(C)$.

\medskip

 (3). Apply the family switching formula [Liu3] 
to the local family invariant contribution 
 and identify it with certain mixed family invariant of $C-\sum_{1\leq i\leq k}
 n_ie_i$ over $\overline{S_{\cal C}}$. The changing of the 
multiplicities of $e_i$, $(1, 1, 1, \cdots, 1)$ to 
 $(n_1, n_2, n_3, \cdots, n_k)$ is called the switching process which
 has been discussed in section \ref{section; point}.

\medskip

(4). As different $\overline{S_{\cal C}}$ may intersect, the naive 
subtraction of all the local family invariant contributions from 
 the original ${\cal AFSW}_{{\cal X}\mapsto B}(1, C)$ leads to over-subtraction.

\medskip

 Inductively, one has to define a version of modified family invariant 
(see section 5.3. of [Liu1] for some explicit examples in the
 differentiable category) for each $S_{\cal C}$.
Schematically it involves subtracting the modified family invariants
(which have already been defined by the induction hypothesis) of 
 $S_{{\cal C}'}$, $\overline{S_{\cal C}}\supset S_{{\cal C}'}$ to define
 the modified family invariant of $S_{\cal C}$.

\medskip

(5). The inductive scheme passes through $B-\coprod_{\cal C}S_{\cal C}$ and
 finally 
one may define a modified family invariant of $B-\coprod_{\cal C}S_{\cal C}$
 by subtracting all the modified family invariants of $S_{\cal C}$ defined.
 In the algebraic category, it involves the usage of residual intersection
 theory in [F] in arguing that the whole family moduli space of $C$ over $B$, 
${\cal M}_C\mapsto B$, can be separated into the subscheme over 
 $\coprod_{{\cal C}}S_{\cal C}$ and the residual portion (which
 definition involves blowing ups inductively). The
modified invariant of $C$ is actually equal to the intersection number
 defined by the residual intersection theory. (See [Liu1] for a 
discussion in the ${\cal C}^{\infty}$ category.)

\medskip

\begin{defin}\label{defin; admissible}
Let $({\cal C}, S_{\cal C})$ be a pair and let $e_i$ be the generators of
 ${\cal C}$. Suppose that the decomposition 
$(C-\sum_{i\leq k}e_i, \sum_{i\leq k}e_i)$ of $C$ satisfies the 
dimension inequality

  $$d_{GT}(C-\sum_{i\leq k}e_i)+\sum_{i\leq k}d_{GT}(e_i)\geq d_{GT}(C),$$
 then we consider all the $k-$tuples $(n_1, n_2, \cdots, n_k)$, $n_i\in {\bf N}$
 such that

$$d_{GT}(C-\sum_{i\leq k}n_ie_i)+\sum_{i\leq k}d_{GT}(e_i)\geq d_{GT}(C).$$

Moreover, 
there exists at least one $i_0$ with $1\leq i_0\leq k$, $n_{i_0}\geq 2$ such that
 
$$d_{GT}(C-\sum_{i\leq k; i\not=i_0}n_ie_i-(n_{i_0}-1)e_{i_0})+
\sum_{i\leq k}d_{GT}(e_i)\geq d_{GT}(C).$$

All such $k-$tuples define decompositions $(C-\sum_{i\leq k}n_ie_i, \sum_{i\leq k}
n_ie_i)$ which can be derived from $(C-\sum_{i\leq k}e_i, \sum_{i\leq k}e_i)$
 through elementary effective moves maintaining their family dimensions above 
the lower bound ${1\over 2}d_{GT}(C)+p_g+dim_{\bf C}B$.
 The set of all such decompositions forms an equivalence class called the
 admissible decomposition class associated to $({\cal C}, S_{\cal C})$.
\end{defin}

\medskip

  The family switching formula [Liu3] 
allows us to relate the family invariants of
 different $(C-\sum_{i\leq k}n_ie_i, \sum_{i\leq k}n_ie_i)$ over
 $\overline{S_{\cal C}}$. That is why we
 view these different decompositions as equivalent.

\medskip

 In the published long paper [Liu1], one takes the universal spaces projection
 $f_n:M_{n+1}\mapsto
 M_n$ of an algebraic surface $M$ as the fiber bundle ${\cal X}\mapsto B$.
 Let $C$ be a $(1, 1)$
class on $M$, then the set $Q$ has been implicitly chosen to be the set
 of all type $I$ exceptional classes $e$ with $e\cdot (C-{\bf M}(E)E)<0$.
 A discussion parallel to the five steps above has been developed 
 in the ${\cal C}^{\infty}$ category.  Please consult [Liu1] for the
 notations and the details.

\bigskip

\begin{rem}\label{rem; nonideal}
 We have made the simplifying assumptions \ref{assum; perfect} to study the
 stratification of the base space $B$. If the conditions (ii)-(v) in
 assumption \ref{assum; perfect} are not satisfied, the existence loci
 of exceptional curves may not be of the right dimension and one has to
 work with the Kuranishi models of the exceptional class to construct 
the algebraic cycle class representing the moduli cycles.

 Moreover, the violation of the conditions (vi)-(vii) indicates that 
 the moduli spaces of different exceptional 
classes may not intersect transversally in $B$. In the 
algebraic category we may use the intersection cycle class
 of the moduli cycle classes to represent the co-existence of different
 exceptional classes. The theory will be developed else where.
\end{rem}

{}

\end{document}